\documentclass[10pt,reqno]{article}
\pdfoutput=1
\usepackage{geometry}
\usepackage{graphicx}

\usepackage{tocloft}

\usepackage{amssymb}
\usepackage{amsmath}
\usepackage{amsthm}

\usepackage{tikz}

\usepackage{float}
\usepackage{caption}
\usepackage{scalerel}
\usepackage{enumitem}
\usepackage{hyperref}
\hypersetup{
  colorlinks   = true,
  urlcolor     = blue,
  linkcolor    = blue,
  citecolor   = red
}

\numberwithin{equation}{section}
\allowdisplaybreaks

\newtheorem{theorem}{Theorem}[section]
\newtheorem{definition}[theorem]{Definition}
\newtheorem{remark}[theorem]{Remark}
\newtheorem{lemma}[theorem]{Lemma}
\newtheorem{proposition}[theorem]{Proposition}

\newcommand{\circg}{\mathring{g} \makebox[0ex]{}}
\newcommand{\uu}{\underline{u} \makebox[0ex]{}}
\newcommand{\uL}{\underline{L} \makebox[0ex]{}}
\newcommand{\uchi}{\underline{\chi} \makebox[0ex]{}}
\newcommand{\hatchi}{\hat{\chi} \makebox[0ex]{}}
\newcommand{\hatuchi}{\hat{\underline{\chi}} \makebox[0ex]{}}
\newcommand{\tr}{\mathrm{tr} \makebox[0ex]{}}
\newcommand{\R}{\mathrm{R} \makebox[0ex]{}}
\newcommand{\ed}{\mathrm{d} \makebox[0ex]{}}
\newcommand{\diver}{\mathrm{div} \makebox[0ex]{}}
\newcommand{\barv}{\bar{v} \makebox[0ex]{}}
\newcommand{\dvol}{\mathrm{dvol} \makebox[0ex]{}}
\newcommand{\calL}{\mathcal{L} \makebox[0ex]{}}
\newcommand{\circnabla}{\mathring{\nabla} \makebox[0ex]{}}
\newcommand{\bdlangle}{\raisebox{-0.15ex}[0ex][0ex]{$\scaleto{\boldsymbol{\langle }}{3ex}$} \makebox[0ex]{}}
\newcommand{\bdrangle}{\raisebox{-0.15ex}[0ex][0ex]{$\scaleto{\boldsymbol{\rangle}}{3ex}$} \makebox[0ex]{}}

\title{\textsc{On effective uniformization of $2$-sphere and the stability}}
\author{Pengyu Le}
\newcommand{\Address}{{
  \bigskip
  \footnotesize
  \textsc{Beijing Institute of Mathematical Sciences and Applications, Beijing, China}
  
  \textit{E-mail address}: \texttt{pengyu.le@bimsa.cn}
}}
\date{}

\begin{document}

\maketitle

\begin{abstract}

We provide a proof of effective uniformization for nearly round $2$-spheres, utilizing an identity related to the third-order differential of the conformal factor. This identity is connected to the geometry of the embedded spacelike surface within the Minkowski lightcone. Additionally, we investigate the stability of the effective uniformization introduced by Klainerman and Szeftel. Our proof is based on a geometric insight: an isometric embedding of a round sphere into Euclidean space can be constructed using an orthogonal basis of the first eigenspace of the Laplacian operator, with the rectangular coordinates corresponding to the basis functions. By adopting these approaches, we simplify both the proofs of effective uniformization and its stability, while also refining the assumptions underlying both results.

\end{abstract}

\tableofcontents

\section{Introduction}\label{sec 1r}

The uniformization theorem for the $2$-sphere asserts that the metric of the sphere can be conformally deformed into a round metric with constant, positive Gauss curvature. The effective uniformization theorem, as introduced in \cite{KS2022}, provides a bound on the conformal factor of this deformation, expressed in terms of the deviation of the Gauss curvature from the constant value of 1.

In their study of the uniformization problem for nearly round spheres, \cite{KS2022} derived effective results for the conformal factor under the assumption that the Gauss curvature is close to 1, either in the $ L^{\infty} $ or $ H^s $ spaces, with $ s \geq 0 $. The proof utilizes the Onofri inequality from \cite{O1982} and its refinement in \cite{CY1987}.  These findings were subsequently applied to construct general covariant modulated spheres in perturbations of Kerr spacetime. Additional effective uniformization results were reported in \cite{SWW2009} (in H\"older space) and \cite{W2023} (in $ L^2 $ space).

Meanwhile, \cite{L2018} (Theorem 2.2.3) presented a similar result but with a stronger assumption on the Gauss curvature, specifically in $ W^{n,p} $ spaces, where $ n \geq 1 $, $ p > 2 $, or $ n \geq 2 $, $ p > 1 $. This result was applied to study the asymptotic geometry of constant mass aspect function foliations at past null infinity in a vacuum-perturbed Schwarzschild spacetime, though no proof was provided.

In Section \ref{sec 2r}, we present a proof of the effective uniformization theorem from \cite{L2018} for nearly round spheres, utilizing a hidden third-order equation (Equation \eqref{eqn 2.1r}) satisfied by the conformal factor. This approach yields an even stronger result, which depends on the $ L^{\infty} $ norm and the $ L^{p>2} $ norm bounds of the Gauss curvature (Section \ref{sec 3r}), and also implies the effective uniformization theorem of \cite{KS2022}.

Equation \eqref{eqn 2.1r} has a natural geometric interpretation related to the geometry of the embedding sphere within the Minkowski lightcone (Section \ref{sec 4r}). This perspective unveils an intriguing connection between the conformal geometry of the sphere and the Minkowski lightcone. Through this connection, we derive the hidden third-order equation \eqref{eqn 2.1r} satisfied by the conformal factor, which aids in its estimation.

Beyond the effective uniformization theorem, \cite{KS2022} also investigated the stability of effective uniformizations for metrics that are close to one another, proving that if two metrics are sufficiently close, there exists a rigid motion such that the isometry between their effective uniformizations is close to this rigid motion. This result builds upon previous work by \cite{FJM2002}, which improved a classical result from \cite{J1961}.

In the second part of our paper (Sections \ref{sec 5r}-\ref{sec 6r}), we provide an alternative proof of the stability of effective uniformization. This proof is based on the observation that the isometric embedding of a round metric with constant Gauss curvature 1 into Euclidean space is given by an orthogonal basis of the first eigenspace of the Laplacian operator, with the rectangular coordinate functions as the basis. Essentially, if two metrics are close, the orthogonal projection between their corresponding Laplacian eigenspaces is close to an isometry. This idea was previously used in \cite{L2018} to study the Laplace equation for scalar functions on a nearly round sphere in a perturbed Schwarzschild spacetime.

Using this framework, we establish the stability of the isometric embedding of the round sphere into Euclidean space, which is analogous to Proposition 4.3 on almost isometries of $ \mathbb{S}^2 $ in \cite{KS2022}. The advantage of this approach lies in its simplicity, as it avoids the need to reference the work of \cite{FJM2002}, although it does require a higher regularity assumption. Furthermore, we extend the stability result of uniformization to the Sobolev space $ W^{2,p} $, benefiting from the simplicity and explicitness of the isometric embedding given by the eigenfunctions. This approach also leads to a coordinate-independent result, which addresses the question raised in \cite{KS2022}.

In the following, we list some notations used in the paper.
\begin{enumerate}[label=\raisebox{0.1ex}{\scriptsize$\bullet$}]
\item $\| f \|_{g,n,p}$: The $W^{n,p}$ Sobolev norm of the function $f$ on a manifold $\Sigma$ with respect to the metric $g$, defined as
\begin{align*}
	\| f \|_{g,n,p}
	=
	\big( \sum_{k=1}^n \int_{\Sigma} \big| \nabla^k f \big|_{g}^p \big)^{\frac{1}{p}},
\end{align*}
where $\nabla$ denotes the covariant derivative on $(\Sigma, g)$.

$\| \omega \|_{g,n,p}$: The $W^{n,p}$ Sobolev norm of a tensor field $\omega$ is defined similarly.

\item
$H_1^n(\mathbb{S}^2, \circg)$, or simply $H_1^n$: The Hilbert space of $1$-forms on the unit sphere $(\mathbb{S}^2, \circg)$.

\item
$\widehat{H}_2^n(\mathbb{S}^2, \circg)$, or simply $\widehat{H}_2^n$: The Hilbert space of trace-free symmetric $2$-tensors on the unit sphere $(\mathbb{S}^2, \circg)$.
\end{enumerate}

\section{Effective uniformization theorem with a small $W^{n,p}$ assumption}\label{sec 2r}

We recall the semi-effective uniformization theorem proved in \cite{CK1993} (lemma 2.3.2).
\begin{theorem}\label{thm 2.1r}
Let $(\Sigma, g)$ be a $2$-sphere with Gauss curvature $K$. There exists a conformal transformation $g = \Omega^2 \circg$ of the round metric $\circg$ of radius $1$. Furthermore, the conformal factor $\Omega$ can be chosen such that the quantities
\begin{align*}
	\Omega_m
	=
	\inf_{\mathbb{S}^2} r \Omega,
	\quad
	\Omega_M
	=
	\sup_{\mathbb{S}^2} r \Omega,
	\quad
	\Omega_1
	=
	\sup_{\mathbb{S}^2} \Omega^2 \vert \nabla \Omega \vert_{g},
	\quad
	\Omega_{2,p}
	=
	( \int_{\mathbb{S}^2} \Omega^{-3p+2} \vert \nabla^2 \Omega \vert^p \dvol_g)^{1/p}
\end{align*}
are bounded, depending only on
\begin{align*}
	\sup r^2 |K|,
	\quad
	\frac{\mathrm{diam}(\Sigma, g)}{r},
\end{align*}
or if $k>0$, only on
\begin{align*}
	k_m
	=
	\min_{\mathbb{S}^2} r^2 K,
	\quad
	k_M
	=
	\max_{\mathbb{S}^2} r^2 K.
\end{align*}
\end{theorem}
The uniformization theorem above provides an estimate for the conformal factor $\Omega$, which serves as a preliminary estimate. This estimate will be refined using the effective uniformization theorem. We now restate the effective uniformization theorem from \cite{L2018} for a nearly round sphere as follows:

\begin{theorem}\label{thm 2.2r}
Let $(\Sigma, g)$ be a 2-sphere with Gauss curvature $K$. Suppose that the Gauss curvature is close to 1 in the Sobolev space $W^{n,p}(\Sigma, g), n\geq n_p$, where $n_p = \big\{\begin{aligned} &1, && p>2 \\  &2, && 2 \geq p >1 \end{aligned}  \big.$, and that
\begin{align*}
	|K - 1| \leq \epsilon,
	\quad
	\Vert K - 1 \Vert^{n,p}_{g} \leq \epsilon.
\end{align*}
Then, there exists a small positive constant $\epsilon_o$ depending on $n, p$, and a constant $c(n, p)$, such that if $\epsilon \leq \epsilon_o$, then for every point $q \in \Sigma$, there exists a unique conformal factor $\Omega^2$ satisfying
\begin{align*}
	\Omega(q) = 1, 
	\quad
	\ed \Omega(q) &= 0,
\end{align*}
such that $\Omega^{-2} g$ is a round metric with constant Gauss curvature $1$, and
\begin{align*}
	|\log \Omega| \leq c(n,p) \epsilon,
	\quad
	\Vert \log \Omega \Vert^{n+2,p}_{g} \leq c(n,p) \epsilon.
\end{align*}
\end{theorem}

The proof of the aforementioned effective uniformization theorem relies on the following lemma, which provides a divergence identity involving the third derivatives of the conformal factor.
\begin{lemma}\label{lem 2.3r}
Let $(\Sigma, \circg)$ denote the standard round sphere of radius 1, and let $(\Sigma, g = \Omega^2 \circg)$ be a conformal deformation of this metric, where the Gauss curvature is given by $K$. Then, the conformal factor $\Omega$ satisfies the equation
\begin{align*}
    K = \Omega^{-2} \big( 1 - \Delta_{\circg} \log \Omega \big).
\end{align*}
Define the symmetric 2-tensor $\Xi$ by
\begin{align*}
    \Xi = - \circg + |\ed \log \Omega|_{\circg}^2 \, \circg - 2 \ed (\log \Omega) \otimes \ed (\log \Omega) + 2 \nabla_{\circg}^2 \log \Omega.
\end{align*}
Decompose $\Xi$ into its trace and trace-free parts with respect to the metric $\circg$:
\begin{align*}
    \Xi = \widehat{\Xi} + \frac{1}{2} \tr_{\circg} \Xi \circg.
\end{align*}
It follows that
\begin{align*}
    K = - \frac{1}{2\Omega^2} \tr_{\circg} \Xi = - \frac{1}{2} \tr_g \Xi,
\end{align*}
and the divergence identity for $\Xi$ is given by
\begin{align}
    \diver_{\circg} \widehat{\Xi} + \Omega^2 \ed K = 0.
    \label{eqn 2.1r}
\end{align}
\end{lemma}
The divergence identity \eqref{eqn 2.1r} can be interpreted as a nonlinear Bochner formula for the conformal factor $\Omega$, and it can be verified directly. In fact, the tensor $\Xi$ admits a geometric interpretation, and the divergence identity \eqref{eqn 2.1r} follows naturally from this interpretation. Therefore, we defer the proof of Lemma \ref{lem 2.3r} to section \ref{sec 4r}. We now proceed to apply the above divergence identity to prove Theorem \ref{thm 2.2r}.

\begin{proof}[Proof of Theorem \ref{thm 2.2r}]
Assume that $\epsilon < \frac{1}{2}$. Let ${\Omega'}^2$ be the conformal factor constructed using the bounds provided in Theorem \ref{thm 2.1r}. Define ${\circg'} = {\Omega'}^{-2} g$ as the round metric of radius $1$. Then, there exists a unique conformal factor $Q$ of $\circg'$ such that $Q^2 \circg'$ is also a round sphere of radius 1, and it satisfies the following conditions:
\begin{align*}
	Q(q) 
	= 
	{\Omega'}(q), 
	\quad
	\ed \log Q(q) 
	= 
	\ed \log {\Omega'}(q).
\end{align*}
Since ${\Omega'}_m$, ${\Omega'}_M$, ${\Omega'}_1$, and ${\Omega'}_{2,p}$ are all bounded by a constant $c$, according to Theorem \ref{thm 2.1r}, it follows that $Q_m$, $Q_M$, $Q_1$, and $Q_{2,p}$ are also bounded by a constant $c$. Here, we abuse the notation $c$ to represent any absolute constant, which may vary depending on the context.

Define $\Omega = Q^{-1} {\Omega'}$ and $\circg = \Omega^{-2} g$. Then, the following properties hold:
\begin{align*}
	\Omega(q)
	= 
	1, 
	\quad
	\ed \Omega(q) 
	= 
	0,
\end{align*}
and $\circg$ is the round metric of radius 1. Moreover, $\Omega_m$, $\Omega_M$, $\Omega_1$, and $\Omega_{2,p}$ are all bounded by a constant $c$. The uniqueness of $\Omega$ follows from the uniqueness of $Q$.

Next, from the equation satisfied by $\Omega$:
\begin{align*}
	K 
	= 
	\Omega^{-2} (1 - \Delta_{\circg} \log \Omega),
\end{align*}
there exists a constant $c(n,p)$ such that
\begin{align*}
	\Vert \log \Omega \Vert_{\circg,n+2,p} 
	\leq 
	c(n,p).
\end{align*}
In the following, we also abuse the notation $c(n,p)$, which may vary throughout the proof but emphasizes its dependence on $n$ and $p$.

Applying Lemma \ref{lem 2.3r} to $(\Sigma, g)$ and $\Omega$, we obtain
\begin{align*}
	\diver_{\circg} \widehat{\Xi} 
	=
	- \Omega^2 \ed K.
\end{align*}
By the estimates for $\Omega$ and the assumption on $K$, there exists a constant $c(n,p)$ depending on $n$ and $p$ such that
\begin{align*}
	\Vert \Omega^2 \ed K \Vert_{\circg,n-1,p} 
	\leq 
	c(n,p) \epsilon.
\end{align*}
By the elliptic estimate on  $2$-sphere and the Sobolev inequality, we have that
\begin{align*}
	\vert \widehat{\Xi} \vert_{\circg}
	\leq
	c(n,p) \Vert \widehat{\Xi} \Vert_{\circg,n,p}
	\leq
	c(n,p) \epsilon.
\end{align*}

From the definition of $\Xi$, we can obtain an estimate for the Hessian $\nabla_{\circg}^2 \log \Omega$ of $\log \Omega$ in principle. This estimate can then be integrated along geodesics to yield estimates for $\ed \Omega$ and $\Omega$. Recall the expression for $\Xi$:
\begin{align*}
	\Xi 
	= 
	\widehat{\Xi} - \Omega^2 K \circg
	=
	- \circg 
	+ |\ed \log \Omega|_{\circg}^2 \, \circg 
	- 2 \ed (\log \Omega) \otimes \ed (\log \Omega) 
	+ 2 \nabla_{\circg}^2 \log \Omega.
\end{align*}
Let $\gamma(s)$ be a unit-speed geodesic emanating from $q$ in $(\Sigma, \circg)$. Denote by $e_1 = \dot{\gamma}$ the tangent vector and let $e_2$ be the unit normal vector field along $\gamma$ with respect to the metric $\circg$. We now consider the equations derived from $\Xi(e_1, e_2)$ and $\Xi(e_1, e_1)$.

\begin{enumerate}[label=\roman*.]
\item \textit{Equation derived from $\Xi(e_1, e_2)$}. Along the geodesic $\gamma$, we obtain the equation
\begin{align*}
	\frac{\ed}{\ed s} (e_2 \log \Omega) 
	- (e_1 \log \Omega) \cdot e_2 \log \Omega 
	= 
	\frac{1}{2} \widehat{\Xi}(e_1, e_2).
\end{align*}
At $s = 0$, since $\ed \log \Omega = 0$, we have the initial condition $e_2 \log \Omega|_{s=0} = 0$. Solving this equation gives
\begin{align}
	e_2 \log \Omega
	=
	\int_0^s 
	\exp \big[ \int_{s'}^s (e_1 \log \Omega) \, \ed s'' \big] 
	\cdot 
	\frac{1}{2} \widehat{\Xi}(e_1, e_2) \, \ed s'.
\label{eqn 2.2r}
\end{align}
By the bounds on $\ed \log \Omega$ and $\widehat{\Xi}$, we find that
\begin{align*}
	|e_2 \log \Omega|
	\leq 
	c(n,p) \epsilon.
\end{align*}

\item \textit{Equation derived from $\Xi(e_1, e_1)$}. Similarly, we have the equation
\begin{align}
	2 \frac{\ed^2}{\ed s^2} \log \Omega 
	- \big( \frac{\ed}{\ed s} \log \Omega \big)^2 
	+ \Omega^2 - 1 
	= 
	-(e_2 \log \Omega)^2 
	+ \widehat{\Xi}(e_1, e_1) 
	+ (1-K) \Omega^2,
\label{eqn 2.3r}
\end{align}
with the initial conditions $\log \Omega|_{s=0} = 0$ and $\frac{\ed}{\ed s} \log \Omega|_{s=0} = 0$. The right-hand side of this equation is bounded as
\begin{align*}
	\vert 
		-(e_2 \log \Omega)^2 
		+ \widehat{\Xi}(e_1, e_1) 
		+ (1-K) \Omega^2 
	\vert
	\leq
	c(n,p) \epsilon.
\end{align*}
Since both $\log \Omega$ and $\frac{\ed}{\ed s} \log \Omega$ are bounded, we can apply the dependence of the solution to the ordinary differential equation on the coefficients to conclude that for sufficiently small $\epsilon$, there exists a constant $c(n,p)$ such that
\begin{align*}
	\vert \log \Omega \vert
	\leq 
	c(n,p) \epsilon, 
	\quad
	\big\vert \frac{\ed}{\ed s} \log \Omega\big\vert
	\leq
	c(n,p) \epsilon.
\end{align*}
\end{enumerate}
The remainder of the theorem follows from the equation $K = \Omega^{-2}(1 - \Delta_{\circg} \log \Omega)$, the assumption on $K$, and the above bound on $\Omega$.
\end{proof}

\begin{remark}\label{rem 2.4r}
We observe that, in order to obtain the $L^\infty$ bound for $\Xi$ via the Sobolev inequality, the regularity assumption on the Gauss curvature should be chosen as in Theorem \ref{thm 2.2r}. However, note that to estimate $\log \Omega$ through equations \eqref{eqn 2.2r} and \eqref{eqn 2.3r}, the $L^\infty$ norm of $\widehat{\Xi}$ is more than what is strictly necessary. Heuristically, to control $|\log \Omega|$ via Sobolev embedding, we need its Sobolev norm $\|\log \Omega\|_{\circg,2,p}$ for some $p > 1$, where the leading term is roughly given by $\|\Xi\|_{\circg,0,p}$, which can be derived from the $L^\infty$ norm of the Gauss curvature $K$. Thus, this approach can be employed to prove the effective uniformization result in \cite{KS2022}. We will discuss this in more detail in the next section.
\end{remark}

We present the following effective uniformization result, which only requires that $K-1$ is small in the $L^\infty$ norm, in a manner comparable to the effective uniformization result in \cite{KS2022}.
\begin{theorem}\label{thm 2.5r}
Let $(\Sigma, g)$ be a 2-sphere with Gauss curvature $K$. Suppose that
\begin{align*}
	|K - 1| \leq \epsilon.
\end{align*}
Then, there exists a positive constant $\epsilon_o$ and a constant $c$, such that if $\epsilon \leq \epsilon_o$, then for every point $q \in \Sigma$, there exists a unique conformal factor $\Omega^2$ satisfying
\begin{align*}
	\Omega(q) = 1, 
	\quad
	\ed \Omega(q) &= 0,
\end{align*}
such that $\Omega^{-2} g$ is a round metric with constant Gauss curvature $1$, and
\begin{align*}
	|\Omega -1| \leq c \epsilon.
\end{align*}
\end{theorem}

In the next section, we will prove a stronger effective uniformization result that requires a bounded 
$L^{\infty}$ norm on $K$. Theorem \ref{thm 2.5r} will then follow as a corollary.

\section{Effective uniformization theorem with bounded $L^{\infty}$ and $L^{p>2}$ assumptions}\label{sec 3r}

We present the following effective uniformization theorem, which assumes that $K$ is bounded.

\begin{theorem}\label{thm 3.1r}
Let $p>2$. Let $(\Sigma, g)$ be a 2-sphere with area $4\pi$ and Gauss curvature $K$. Then, there exists a constant $c_p$, not only on $p$, but also on
\begin{align*}
	\sup |K|,
	\quad
	\mathrm{diam}(\Sigma, g),
\end{align*}
or if $k>0$, on
\begin{align*}
	k_m
	=
	\min_{\mathbb{S}^2} K,
	\quad
	k_M
	=
	\max_{\mathbb{S}^2} K,
\end{align*}
such that for every point $q \in \Sigma$, there exists a unique conformal factor $\Omega^2$ satisfying
\begin{align*}
	\Omega(q) = 1, 
	\quad
	\ed \Omega(q) &= 0,
\end{align*}
where $\Omega^{-2} g$ is a round metric with constant Gauss curvature $1$, and
\begin{align*}
	\| \Omega -1 \|_{g,2,p} \leq c_p \| K - 1 \|_{g,0,p}.
\end{align*}
\end{theorem}

It is evident that theorem \ref{thm 2.5r} follows as a corollary of theorem \ref{thm 3.1r}.
The key to the proof is to first obtain the $L^{p}$ estimate for $\widehat{\Xi}$ using the equation $\diver_{\circg} \widehat{\Xi} = - \Omega^2 \ed K$, given the $L^{p}$ norm of $K-1$. The following lemma provides the desired estimate.

\begin{lemma}\label{lem 3.2r}
Suppose $\widehat{\Xi}$ satisfies the equation
\begin{align*}
	\diver_{\circg} \widehat{\Xi} = - \Omega^2 \ed K,
\end{align*}
\begin{enumerate}[label=\roman*.]
\item {\bf $L^2$ estimate:} There exists a constant $c_{p,q}$ such that
\begin{align*}
	\| \widehat{\Xi} \|_{\circg,0,2}
	\leq
	c_{p,q} \| \Omega^2 \|_{\circg,1, q} \cdot \| K - k \|_{\circg,0,p},
	\quad
	\forall k \in \mathbb{R}.
\end{align*}
where $\frac{1}{p}+\frac{1}{q} <1$, $p\geq 2$.

\item {\bf $L^{p>2}$ estimate:} There exists a constant $c_{q,r}$ such that
\begin{align*}
	\| \widehat{\Xi} \|_{\circg,0,p}
	\leq
	c_{q,r} \| \Omega^2 \|_{\circg,1, r} \cdot \| K - k \|_{\circg,0,q},
	\quad
	\forall k \in \mathbb{R},
\end{align*}
where $\frac{1}{q} + \frac{1}{r} = \frac{1}{2} + \frac{1}{p}$, $q \geq p$.

\end{enumerate}
\end{lemma}

\begin{proof}
Define the operator $\calL$ from $H_1^1$ to $\widehat{H}_2^0$ as
\begin{align*}
	(\calL \omega)_{ij}
	=
	\frac{1}{2} ( \circnabla_i \omega_j + \circnabla_j \omega_i) 
	- \frac{1}{2} \diver_{\circg} \omega \cdot \circg.
\end{align*}
The operator $\calL$ is the $L^2$-adjoint of $\diver_{\circg}$. By the divergence theorem, we have
\begin{align*}
	\int_{\Sigma} \langle \widehat{\Xi}, \calL \omega \rangle_{\circg} \dvol_{\circg}
	=
	\int_{\Sigma} \langle \diver_{\circg} \widehat{\Xi}, \omega \rangle_{\circg} \dvol_{\circg}
	=
	\int_{\Sigma} \langle - \Omega^2 \ed K, \omega \rangle_{\circg} \dvol_{\circg}
	=
	\int_{\Sigma} \langle K - k,  \diver_{\circg} (\Omega^2 \omega) \rangle_{\circg} \dvol_{\circg}
\end{align*}
\begin{enumerate}[label=\textit{\roman*}.]
\item {\bf $L^2$ estimate:} 
We obtain the following estimate:
\begin{align*}
	&\phantom{\leq}
	\big| \int_{\Sigma} \langle \widehat{\Xi}, \calL \omega \rangle_{\circg} \dvol_{\circg} \big|
\\	
	&
	\leq
	\big| \int_{\Sigma} (K - k) \Omega^2 \diver_{\circg} \omega \rangle_{\circg} \dvol_{\circg} \big|
	+
	\big| \int_{\Sigma} (K - k) \langle  \ed (\Omega^2), \omega \rangle_{\circg} \dvol_{\circg} \big|
\\
	&
	\leq
	c \| K - k \|_{\circg,0,p} \cdot \| \Omega^2 \|_{\circg,0,q'} \cdot \|  \ed \omega \|_{\circg,0,2}
	+
	c \| K - k \|_{\circg,0,p} \cdot \| \ed (\Omega^2) \|_{\circg,1,q} \cdot \|  \omega \|_{\circg,0, r}
\\
	&
	\leq
	c_{p,q} \| K - k \|_{\circg,0,p} \cdot \| \Omega^2 \|_{\circg,1,q} \cdot \| \omega \|_{\circg,1,2}
\end{align*}
where $\frac{1}{p} + \frac{1}{q'} + \frac{1}{2} < 1$, $\frac{1}{q'} = \frac{1}{q} -\frac{1}{2}$, $\frac{1}{p} + \frac{1}{q} + \frac{1}{r} = 1$.

Note that $\calL$ is surjective since its $L^2$-adjoint, $\diver_{\circg}$, maps from $H_2^0$ to $H_1^{-1}$ and is injective by the elliptic estimate for $\diver_{\circg}$. Therefore, we can write
\begin{align*}
	\| \widehat{\Xi} \|_{\circg,0,2}
	=
	\sup_{\calL \omega \neq 0}
	\frac{\langle \widehat{\Xi}, \calL \omega \rangle_0}{ \| \calL \omega \|_{\circg,0,2} }.
\end{align*}
Next, introduce the closed linear subspace $L^1 \subset H_1^1$, spanned by the set $\{ \ed Y_l^m \}_{l \in \mathbb{Z}_{\geq 2}, -l \leq m \leq l } \cup \{ \star\ed Y_l^m \}_{l \in \mathbb{Z}_{+}, -l \leq m \leq l }$. Then, the operator $\calL: L^1 \to \widehat{H}_2^0$ is a self-adjoint bijection, and there exists a constant $c$ such that
\begin{align*}
	\sup_{\omega \neq 0, \omega \in L^1} \frac{\| \calL \omega \|_{\circg,0,2}}{\| \omega \|_{\circg,1,2}},
	\quad
	\sup_{\omega \neq 0, \omega \in L^1} \frac{\| \omega \|_{\circg,1,2}}{\| \calL \omega \|_{\circg,0,2}}
	\leq
	c.
\end{align*}
Therefore, we conclude that
\begin{align*}
	\| \widehat{\Xi} \|_{\circg,0,2}
	&=
	\sup_{\calL \omega \neq 0}
	\frac{\langle \widehat{\Xi}, \calL \omega \rangle_0}{ \| \calL \omega \|_{\circg,0,2} }
	\leq
	c_{p,q} \| K - k \|_{\circg,0,p} \cdot \| \Omega^2 \|_{\circg,1,q}
	\cdot 
	\sup_{\omega \neq 0, \omega \in L^1}
	\frac{ \| \omega \|_{\circg,1,2}}{ \| \calL \omega \|_{\circg,0,2} }
\\
	&
	\leq
	c_{p,q} \| K - k \|_{\circg,0,p} \cdot \| \Omega^2 \|_{\circg,1,q}.
\end{align*}

\item 
{\bf $L^{p>2}$ estimate:} Similarly,
\begin{align*}
	&\phantom{\leq}
	\big| \int_{\Sigma} \langle \widehat{\Xi}, \calL \omega \rangle_{\circg} \dvol_{\circg} \big|
\\	
	&
	\leq
	\big| \int_{\Sigma} (K - k) \Omega^2 \diver_{\circg} \omega \rangle_{\circg} \dvol_{\circg} \big|
	+
	\big| \int_{\Sigma} (K - k) \langle  \ed (\Omega^2), \omega \rangle_{\circg} \dvol_{\circg} \big|
\\
	&
	\leq
	c \| K - k \|_{\circg,0,q} \cdot \| \Omega^2 \|_{\circg,0,r'} \cdot \|  \ed \omega \|_{\circg,0,p'}
	+
	c \| K - k \|_{\circg,0,q} \cdot \| \ed (\Omega^2) \|_{\circg,1,r} \cdot \|  \omega \|_{\circg,0,p''}
\\
	&
	\leq
	c_{p,q} \| K - k \|_{\circg,0,q} \cdot \| \Omega^2 \|_{\circg,1,r} \cdot \| \omega \|_{\circg,1,p'}
\end{align*}
where $\frac{1}{p} + \frac{1}{p'} = 1$, $\frac{1}{r'}= \frac{1}{r} - \frac{1}{2}$, $\frac{1}{p''} = \frac{1}{p'} - \frac{1}{2}$.

From $p'<2$ and $\calL$ being surjective from $H_1^1$ to $\widehat{H}_2^0$, we obtain that $\calL$ is surjective from $W_1^{1,p'}$ to $\widehat{W}_2^{0,p'}$. Therefore
\begin{align*}
	\| \widehat{\Xi} \|_{\circg,0,p}
	=
	\sup_{\calL \omega \neq 0}
	\frac{\langle \widehat{\Xi}, \calL \omega \rangle_0}{ \| \calL \omega \|_{\circg,0,p'} }.
\end{align*}
Similarly introduce the closed linear subspace $L^{1,p'} \subset W_1^{1,p'}$, spanned by the set $$\{ \ed Y_l^m \}_{l \in \mathbb{Z}_{\geq 2}, -l \leq m \leq l } \cup \{ \star\ed Y_l^m \}_{l \in \mathbb{Z}_{+}, -l \leq m \leq l }.$$
Then, the operator $\calL: L^{1,p'} \rightarrow \widehat{W}_2^{0,p'}$ is a continuous bijection, and there exists a constant $c_{p}$ such that
\begin{align*}
	\sup_{\omega \neq 0, \omega \in L^{1,p'}} \frac{\| \calL \omega \|_{\circg,0,p'}}{\| \omega \|_{\circg,1,p'}},
	\quad
	\sup_{\omega \neq 0, \omega \in L^{1,p'}} \frac{\| \omega \|_{\circg,1,p'}}{\| \calL \omega \|_{\circg,0,p'}}
	\leq
	c_p.
\end{align*}
Therefore, we conclude that
\begin{align*}
	\| \widehat{\Xi} \|_{\circg,0,p}
	&=
	\sup_{\calL \omega \neq 0}
	\frac{\langle \widehat{\Xi}, \calL \omega \rangle_0}{ \| \calL \omega \|_{\circg,0,p'} }
	\leq
	c_{q,r} \| K - k \|_{\circg,0,q} \cdot \| \Omega^2 \|_{\circg,1,r}
	\cdot 
	\sup_{\omega \neq 0, \omega \in L^{1,p'}}
	\frac{ \| \omega \|_{\circg,1,p'}}{ \| \calL \omega \|_{\circg,0,p'} }
\\
	&
	\leq
	c_{q,r} \| K - k \|_{\circg,0,q} \cdot \| \Omega^2 \|_{\circg,1,r}.
\end{align*}
\end{enumerate}
This completes the proof of the lemma.
\end{proof}

Now we turn to the proof of theorem \ref{thm 3.1r}.
\begin{proof}[Proof of theorem \ref{thm 3.1r}]
Analogous to the proof of Theorem \ref{thm 2.2r}, there exists a unique conformal factor $\Omega$ such that $\circg = \Omega^{-2} g$ is a round metric of radius 1, with $\Omega(q) = 1$ and $\ed \Omega(q) = 0$. Furthermore, the quantities $\Omega_m$, $\Omega_M$, $\Omega_1$, and $\Omega_{2,p}$ are all bounded by a constant $c$ whose dependence is specified as in the theorem.

By Lemma \ref{lem 3.2r}, we obtain the estimate
\begin{align*}
    \| \widehat{\Xi} \|_{\circg,0,p}
    \leq
    c_p \| K - 1 \|_{\circg,0,p}
    \leq 
    c_p \| K - 1 \|_{g,0,p}.
\end{align*}
To utilize this estimate of $\widehat{\Xi}$, we transform the formula for $\Xi$ into another form. Note that $\log \Omega = - \log \Omega^{-1}$, so we have
\begin{align*}
	2 \circnabla^2 \log \Omega
	=
	-2 \circnabla^2 \log \Omega^{-1}
	=
	-2 \circnabla \big( \frac{\circnabla \Omega^{-1}}{\Omega^{-1}} \big)
	=
	- 2 \frac{\circnabla^2 \Omega^{-1}}{\Omega^{-1}}
	+ 2 \circnabla (\log \Omega) \otimes \circnabla (\log \Omega),
\end{align*}
and therefore
\begin{align*}
	&
	\Xi
	=
	\big[ \vert \ed (\log \Omega^{-1}) \vert^2 - 1 \big] \circg
	- 2 \Omega \circnabla^2 \Omega^{-1}
	=
	-\Omega^2 K \circg
	+ \widehat{\Xi},
\\
	&
	\widehat{\Xi}
	=
	- \Omega \big[ 2 \circnabla^2 \Omega^{-1} - \Delta_{\circg} \Omega^{-1} \cdot \circg \big].
\end{align*}
Following the construction in the proof of Theorem \ref{thm 2.2r}, we now consider the equations derived from $\Xi(e_1, e_2)$ and $\Xi(e_1, e_1)$.

\begin{enumerate}[label=\roman*.]
\item \textit{Equation derived from $\Xi(e_1, e_2)$}. 
Along the geodesic $\gamma$, we obtain the equation
\begin{align*}
	\frac{\ed}{\ed s} (e_2 \Omega^{-1})
	=
	- \frac{1}{2} \Omega^{-1} \widehat{\Xi}(e_1, e_2).
\end{align*}
Next, consider the $L^{2}$ norm of $e_2 \Omega^{-1}$:
\begin{align*}
	&
	\big\vert \frac{\ed}{\ed s} \int_0^{2\pi} \vert e_2 \Omega^{-1} \vert^{2} \ed \phi \big\vert
	\leq
	\big\vert \int_0^{2\pi} \vert e_2 \Omega^{-1} \vert \cdot  \Omega^{-1} \widehat{\Xi}(e_1, e_2) \ed \phi \big\vert
\\
	&\phantom{\vert \frac{\ed}{\ed s} \int_0^{2\pi} \vert e_2 \Omega^{-1} \vert^{2} \ed \phi \vert}
	\leq
	c_p (\sin s)^{-\frac{1}{p}}
	\cdot
	\big( \int_0^{2\pi} \vert e_2 \Omega^{-1} \vert^{2} \ed \phi \big)^{\frac{1}{2}}
	\cdot
	\big( \int_0^{2\pi} \sin s \cdot \vert \widehat{\Xi} \vert^4 \ed \phi \big)^{\frac{1}{p}},
\end{align*}
Therefore, we have the estimate
\begin{align*}
	\big\vert \frac{\ed}{\ed s} \| e_2 \Omega^{-1}(s, \cdot) \|_{\mathbb{S}^1,0,2} \big\vert
	\leq
	c (\sin s)^{-\frac{1}{p}}
	\big( \int_0^{2\pi} \sin s \cdot \vert \widehat{\Xi} \vert^p \ed \phi \big)^{\frac{1}{p}},
\end{align*}
which implies the following estimate by integrating the inequality:
\begin{align*}
	&
	\| e_2 \Omega^{-1}(s, \cdot) \|_{\mathbb{S}^1,2}
	\leq
	c_p \big[ 
		\int_{0}^s (\sin s')^{-\frac{1}{p-1}} \ed s' 
		\big]^{\frac{p-1}{p}}
	\cdot
	\big[ 
		\int_{0}^s  \int_0^{2\pi} \sin s' \cdot \vert \widehat{\Xi} \vert^p \ed \phi  \ed s' 
		\big]^{\frac{1}{p}}.
\\
	\Rightarrow\quad
	&
	\| e_2 \Omega (s, \cdot) \|_{\circg,0,2}
	\leq
	c_p \| K - 1 \|_{g,0,p}.
\end{align*}

\item \textit{Equation derived from $\Xi(e_1, e_1)$}. 
We transform equation \eqref{eqn 2.3r} into an equation for $\Omega^{-\frac{1}{2}}$. Note that
\begin{align*}
	2 \frac{\ed^2}{\ed s^2} \log \Omega 
	- \big( \frac{\ed}{\ed s} \log \Omega \big)^2 
	=
	-4 \frac{\frac{\ed^2}{\ed s^2} \Omega^{-\frac{1}{2}} }{\Omega^{-\frac{1}{2}}}.
\end{align*}
Thus, we have
\begin{align*}
	4 \frac{\ed^2}{\ed s^2} \Omega^{-\frac{1}{2}} 
	- \big( \Omega^{-\frac{1}{2}} \big)^{-3} + \Omega^{-\frac{1}{2}} 
	= 
	\Omega^{-\frac{1}{2}} (e_2 \log \Omega)^2 
	- \Omega^{-\frac{1}{2}} \widehat{\Xi}(e_1, e_1) 
	+ (K-1) \Omega^{\frac{3}{2}}.
\end{align*}
This implies that
\begin{align*}
	&\phantom{=}
	\frac{\ed}{\ed s} 
	\big[ 
		4 \big(  \frac{\ed}{\ed s} \Omega^{-\frac{1}{2}} \big)^2 
		+ \Omega^{-1} (\Omega - 1)^2
	\big]
\\
	&
	=
	2 \frac{\ed}{\ed s} \Omega^{-\frac{1}{2}}
	\cdot
	\big[
		\Omega^{-\frac{1}{2}} (e_2 \log \Omega)^2 
		- \Omega^{-\frac{1}{2}} \widehat{\Xi}(e_1, e_1) 
		+ (K-1) \Omega^{\frac{3}{2}}
	\big].
\end{align*}
Therefore, by integrating the above equation, we obtain
\begin{align*}
	&\phantom{=}
	\big\vert
		\frac{\ed}{\ed s}
		\big\{
			\int_0^{2\pi}
			\big[
				4 \big(  \frac{\ed}{\ed s} \Omega^{-\frac{1}{2}} \big)^2 
				+ \Omega^{-1} (\Omega - 1)^2
			\big] \ed \phi
		\big\}^{\frac{1}{2}}
	\big\vert
\\
	&
	\leq
	\frac{1}{2}	
	\big[
		\int_0^{2\pi}
			\Omega^{-1} (e_2 \log \Omega)^4
		\ed \phi
	\big]^{\frac{1}{2}}
	+
	c_p
	\big[
		\int_0^{2\pi}
			\Omega^{-1} \vert \widehat{\Xi}(e_1, e_1) \vert^p
		\ed \phi
	\big]^{\frac{1}{p}}
	+
	c^p
	\big[
		\int_0^{2\pi}
			\Omega^{\frac{3p}{2}} (K-1)^p
		\ed \phi
	\big]^{\frac{1}{p}}
\\
	&
	\leq
	c_p \| K -1 \|_{g,0,p}
	+
	c_p (\sin s)^{-\frac{1}{p}}
	\big( \int_0^{2\pi}
		\sin s \cdot \vert \widehat{\Xi}\vert^p
		\ed \phi 
	\big)^{\frac{1}{p}}
	+
	c_p (\sin s)^{-\frac{1}{p}}
	\big( \int_0^{2\pi}
		\sin s \cdot \vert K-1 \vert^p
		\ed \phi 
	\big)^{\frac{1}{p}}.
\end{align*}
This implies that
\begin{align*}
	&
	\big\{
		\int_0^{2\pi}
		\big[
			4 \big(  \frac{\ed}{\ed s} \Omega^{-\frac{1}{2}} \big)^2 
			+ \Omega^{-1} (\Omega - 1)^2
		\big] \ed \phi
	\big\}^{\frac{1}{2}}
	\leq
	c_p \| K -1 \|_{g,0,p}
\\
	\Rightarrow
	\quad
	&
	\| e_1 \Omega \|_{\circg,0,2}
	\leq
	c_p \| K -1 \|_{g,0,p},
	\quad
	\| \Omega - 1 \|_{\circg,0,2} 
	\leq
	c_p \| K -1 \|_{g,0,p}.
\end{align*}
\end{enumerate}
Combining i. and ii., we obtain
\begin{align*}
	\| \Omega -1 \|_{\circg, 1, 2}
	\leq
	c_p \| K -1 \|_{g,0,p}
\end{align*}
Then, the formula
\begin{align*}
	2 \Omega \circnabla^2 \Omega^{-1}
	=
	- \widehat{\Xi}
	+ 
	\vert \ed (\log \Omega^{-1}) \vert^2 \circg
	+ \Omega^2 (K -1) \circg
	+ (\Omega^2 - 1) \circg
\end{align*}
implies that
\begin{align*}
	\| \Omega - 1 \|_{\circg, 2,2}
	\leq
	c_p \| K -1 \|_{g,0,p}
	\quad
	\Rightarrow
	\quad
	\vert \Omega - 1 \vert
	\leq
	c_p \| K -1 \|_{g,0,p}.
\end{align*}
Then the theorem is a consequence of the standard elliptic estimate.
\end{proof}

\section{Interpretation of Lemma \ref{lem 2.3r} in Minkowski lightcone}\label{sec 4r}

Following the notation introduced in Lemma \ref{lem 2.3r}, $(\Sigma, g = \Omega^2 \circg)$ represents a conformal deformation of the round metric of radius $1$ on the sphere. We now embed $(\Sigma, g)$ into the lightcone of Minkowski spacetime $\mathbb{M}^{3,1}$ using the following method.

Let $\{t, r, \vartheta\}$ denote the spatial polar coordinate system of $\mathbb{M}^{3,1}$, where $\vartheta$ represents the coordinates on the sphere $\mathbb{S}^{2}$. Introduce a double null coordinate system $\{\uu, u, \vartheta\}$ by defining
\begin{align*}
	\uu = \frac{t + r}{2}, \quad u = \frac{t - r}{2}.
\end{align*}
The coordinate vector fields $\partial_{\uu}$ and $\partial_u$ are null,
\begin{align*}
	\partial_{\uu} = \partial_t + \partial_r, \quad \partial_u = \partial_t - \partial_r.
\end{align*}
Let $C_0$ be the outgoing lightcone from the origin. The pair $\{\uu, \vartheta\}$ forms a coordinate system on $C_0$, where the degenerate metric on $C_0$ is given by $\uu^2 \circg$, with $\circg$ being the standard metric on the unit sphere. Therefore, $(\Sigma, g = \Omega^2 \circg)$ is isometric to the following surface in $C_0$:
\begin{align*}
	\bar{\Sigma} = \{ (\uu, \vartheta) : \uu = \Omega (\vartheta) \}.
\end{align*}

We now define a particular conjugate null frame on $\bar{\Sigma}$, as well as the corresponding null expansions.

\begin{definition}\label{def 4.1r}
Let $L$ be the position vector field of $\Sigma$, i.e., the vector emanating from the origin to the points in $\Sigma$, so that $L(q) = \vec{oq}$. Define the conjugate null normal vector field $\uL$ along $\bar{\Sigma}$ to be the vector field satisfying the conjugate condition
\begin{align*}
	g(\uL, L) = -2.
\end{align*}
The pair $\{L, \uL\}$ forms a conjugate null frame along $\bar{\Sigma}$. A schematic representation of this setup is provided in Figure \ref{fig 1r}.
\begin{figure}[H]
\centering
\begin{tikzpicture}[scale=0.6]
\draw[fill] (0,0) circle (1pt);
\draw (0,0) node[below right] {$o$}  -- (4,4) (0,0) -- (-4,4);
\draw[thick, dashed] (3,3) to [out=45,in=15] (2,3.5)
to [out=-165,in=-40] (-2,3)
to [out=140,in=135] (-3.5,3.5); 
\draw[dashed] 
(0,0) -- (-2,3);
\draw[thick] (-3.5,3.5) to [out=-45,in=-160] (-1.8,3.2)
to [out=20,in=-170] (-1,3.5) 
to [out=10,in=150] (1,2.5)
to [out=-30,in=-135] (3,3); 
\draw 
(0,0) -- (1,2.5); 
\path[fill=black] (1,2.5) circle(2pt);
\draw[->,thick] (1,2.5) -- (2,5) node[above] {\footnotesize$L$};  
\draw[->,thick] (1,2.5) -- (0.3,4) node[above] {\footnotesize$\uL$}; 
\path[fill=black] (3.12,3.12) circle(2pt);
\draw[->,thick] (3.12,3.12) -- (6,6) node[above] {\footnotesize$L$}; 
\draw[->,thick] (3.12,3.12) -- (3.12-1,3.12+1) node[right] {\footnotesize$\uL$}; 
\path[fill=black] (-3.6,3.6) circle(2pt);
\draw[->,thick] (-3.6,3.6) -- (-7.2,7.2) node[above] {\footnotesize$L$}; 
\draw[->,thick] (-3.6,3.6) -- (-2.8,4.5) node[above] {\footnotesize$\uL$}; 
\end{tikzpicture}
\caption{Conjugate null frame $\{L, \uL\}$.}
\label{fig 1r}
\end{figure}
\end{definition}

\begin{definition}\label{def 4.2r}
Associated with the conjugate null frame $\{L, \uL\}$ on $\bar{\Sigma}$, we define the second fundamental forms $\chi$ and $\uchi$ as follows:
\begin{align*}
	\chi(X,Y) = g(D_X L, Y), \quad \uchi(X,Y) = g(D_X \uL, Y),
\end{align*}
where $X, Y \in T\Sigma$ and $D_X$ denotes the covariant derivative along the surface $\Sigma$. Let $g$ denote the intrinsic metric on $\bar{\Sigma}$. We decompose $\chi$ and $\uchi$ into their trace and trace-free parts:
\begin{align*}
	\chi = \hatchi + \frac{\tr \chi}{n-1} g, \quad \uchi = \hatuchi + \frac{\tr \uchi}{n-1} g,
\end{align*}
where $\tr \chi$ and $\tr \uchi$ are the null expansions, and $\hatchi$ and $\hatuchi$ are the null shears.
\end{definition}

The Gauss and Codazzi equations for the null second fundamental forms $\chi$ and $\uchi$ are as follows. Let $\nabla$ denote the covariant derivative on $(\bar{\Sigma}, g)$. Direct calculations show that $\chi = g$.

\begin{enumerate}[label=,leftmargin=15pt]
\item \textbf{The Gauss equation:} Let $\R_{ijkl}$ be the Riemann curvature tensor\footnote{The convention for the Riemann tensor is $\mathrm{R}_{ijkl} = g(\nabla_{\partial_i} \nabla_{\partial_j} \partial_k - \nabla_{\partial_j} \nabla_{\partial_i} \partial_k , \partial_l)$} of $(\Sigma, g)$. Then, the Gauss equation is
\begin{align*}
	\R_{ijkl} 
	= 
	\frac{1}{2} \big( g_{ik} \uchi_{jl} + \uchi_{ik} g_{jl} - g_{jk} \uchi_{il} - \uchi_{jk} g_{il} \big).
\end{align*}

\item \textbf{The Codazzi equation:}
\begin{align*}
	\nabla_i \uchi_{jk} - \nabla_j \uchi_{ik} = 0.
\end{align*}
\end{enumerate}

From these equations, we obtain the following lemma.

\begin{lemma}\label{lem 4.3r}
The Gauss curvature $K$ of $(\bar{\Sigma}, g)$ satisfies the equation
\begin{align*}
	K = -\frac{1}{2} \tr \uchi.
\end{align*}
Additionally, the null second fundamental form $\uchi$ satisfies the equation
\begin{align*}
	\diver \hatuchi - \frac{1}{2} \nabla \tr \uchi = 0.
\end{align*}
\end{lemma}

Finally, Lemma \ref{lem 2.3r} follows from the above lemma and the formula
\begin{align}
	\uchi = \Xi = - \circg + |\ed \log \Omega|_{\circg}^2 \, \circg - 2 \ed (\log \Omega) \otimes \ed (\log \Omega) + 2 \nabla_{\circg}^2 \log \Omega.
	\label{eqn 4.1r}
\end{align}
One can easily verify the above formula through the following steps.
\begin{enumerate}
\item
Let $(\theta^1, \theta^2)$ denote the coordinate system on the sphere. Let $\partial_i$, for $i=1,2$, be the coordinate vector fields corresponding to the coordinate system $(t, r, \theta^1, \theta^2)$. Let $\bar{\partial}_i$ be the coordinate vector fields corresponding to the coordinate system $(\theta^1, \theta^2)$ on $\bar{\Sigma}$. Then, we have the relation
\begin{align*}
	\bar{\partial}_i
	=
	\partial_i
    	+
	(\circnabla_i \Omega) \partial_{\uu}.
\end{align*}

\item
The vector fields $L$ and $\uL$ are given by the following equations:
\begin{align*}
	L
	= 
	\Omega \partial_{\uu}, 
	\quad
	\uL
	= 
	\Omega^{-1} 
	\big[ 
		\partial_u 
		+ \Omega^{-2} \big| \ed \Omega \big|_{\circg}^2 \partial_{\uu} 
		+ 2 \Omega^{-2} (\circg^{-1})^{ij} (\circnabla_j \Omega) \partial_i
	\big].
\end{align*}

\item
Note that the following relations hold:
\begin{align*}
	g(D_{\partial_i} \partial_j, \partial_k) 
	= 
	\Omega^2 \circg (\circnabla_{\partial_i} \partial_j, \partial_k),
	\quad
	D_{\partial_i} \partial_{\uu} 
	= 
	- D_{\partial_i} \partial_u 
	= 
	\partial_i.
\end{align*}
Equation \eqref{eqn 4.1r} follows directly from the definition of $\uchi$.
\end{enumerate}

\section{Stability of isometric embedding of round sphere to Euclidean space}\label{sec 5r}

In this section, we examine a special case of the stability of effective uniformizations, where both metrics are round with Gauss curvature equal to $1$. The general case will follow from the effective uniformization theorem, building on this special case.

Let $g_{1}$ and $g_{2}$ denote two round metrics with Gauss curvature $1$ on the sphere $\Sigma$. There exist isometric embeddings of $(\Sigma, g_{1})$ and $(\Sigma, g_{2})$ into $(\mathbb{S}^2, \circg)$, the round sphere of radius $1$ centered at the origin in the three-dimensional Euclidean space $\mathbb{E}^3$. Let $\Phi_a$ represent the embedding of $(\Sigma, g_{a})$ into $(\mathbb{S}^2, \circg)$, for $a = 1, 2$. The composition map $\Phi_2 \circ \Phi_1^{-1}$ is then a map on $(\mathbb{S}^2, \circg)$. 

We are interested in the distance between the map $\Phi_2 \circ \Phi_1^{-1}$ and the rigid motions of the round sphere $(\mathbb{S}^2, \circg)$. The result is that, if the metrics $g_{a}$ are sufficiently close, there exists a rigid motion $O$ of $(\mathbb{S}^2, \circg)$ such that $\Phi_2 \circ \Phi_1^{-1}$ is close to $O$.

Both embeddings of $(\Sigma, g_{a})$ for $a = 1, 2$ into $\mathbb{E}^3$ can be described by the three eigenfunctions in an orthogonal basis of the first eigenspace $\Lambda_{k=1}^a$ of the Laplacian operator $\Delta_{a}$ on $(\Sigma, g_{a})$, for $a = 1, 2$. We derive an effective estimate for the orthogonal projection map from the eigenspace $\Lambda_{k=1}^1$ to the eigenspace $\Lambda_{k=1}^2$ with respect to the metric $g_{2}$.

\begin{lemma}\label{lem 5.1r}
Suppose that $\Vert g_{2} - g_{1} \Vert_{g_{1},n,p} \leq \delta <1$ where $n\geq n_p$. Let $v\in \Lambda_{k=1}^1$ and $\bar{v}$ denote the orthogonal projection of $v$ to $\Lambda_{k=1}^2$ with respect to the metric $g_{2}$. Then, there exists a constant $c(n,p)$ such that the following inequalities hold
\begin{align*}
	\Vert v- \barv \Vert^{n+1,p}_{g_{1}}
	\leq
	c(n,p) \delta \Vert v \Vert_{g_{1},n+1,p},
	\quad
	\Vert v- \barv \Vert^{n+1,p}_{g_{2}}
	\leq
	c(n,p) \delta \Vert v \Vert_{g_{2},n+1,p}.
\end{align*}
\end{lemma}
\begin{proof}
Consider the eigenvalue equations for $ \Delta_{a} $:
\begin{align*}
	\Delta_{1} v = 2 v,
	\quad
	\Delta_{2} \bar{v} = 2 \bar{v}.
\end{align*}
Taking their difference, we obtain:
\begin{align*}
	\Delta_{1} (v-\barv) + (\Delta_{1} - \Delta_{2}) \barv = 2 (v-\barv),
	\quad
	\Delta_{2} (v-\barv) + (\Delta_{1} - \Delta_{2}) v = 2 (v-\barv)
\end{align*}
This leads to
\begin{align*}
	(\Delta_{2} - 2) (v-\barv) = (\Delta_{2} - \Delta_{1}) v
\end{align*}
Since $v - \bar{v} \in (\Lambda_{k=1}^2)^{\perp, g_{2}}$, the orthogonal complement of $ \Lambda_{k=1}^2$ with respect to the metric $g_{2}$, we have:
\begin{align*}
	\Vert v - \bar{v} \Vert_{g_{2},n+1,p} 
	\leq 
	c(n,p) \Vert (\Delta_{2} - \Delta_{1}) v \Vert_{g_{2},n-1,p}.
\end{align*}

Let $\nabla_{a}$ denote the covariant derivative of $(\Sigma, g_{a})$, and let $(\Gamma_{a})_{ij}^k$ be the corresponding Christoffel symbols. Define $\triangle_{ij}^k = (\Gamma_{2})_{ij}^k - (\Gamma_{1})_{ij}^k$. Then, we have:
\begin{align*}
	\triangle_{ij}^k 
	= 
	\frac{1}{2} 
	((g_{2})^{-1})^{kl} 
	((\nabla_{1})_i (g_{2})_{jl} + (\nabla_{1})_j (g_{2})_{il} - (\nabla_{1})_l (g_{2})_{ij}),
\end{align*}
and
\begin{align*}
	(\nabla_{2})^2_{ij} v - (\nabla_{1})^2_{ij} v = - \triangle_{ij}^k v_k,
	\quad
	(\Delta_{2} - \Delta_{1}) v 
	= 
	((g_{2})^{-1} - (g_{1})^{-1})^{ij} \cdot (\nabla_{1})^2_{ij} v
	- 
	((g_{2})^{-1})^{ij} \triangle_{ij}^k v_k.
\end{align*}

Now, suppose that $\Vert g_{2} - g_{1} \Vert_{g_{1},n,p} \leq 1$, where $n \geq n_p$. Then, there exists a constant $c(n,p)$ such that:
\begin{align*}
	c(n,p)^{-1} \Vert f \Vert_{g_{1},n+1,p} 
	\leq 
	\Vert f \Vert_{g_{2},n+1,p} 
	\leq 
	c(n,p)\Vert f \Vert_{g_{1},n+1,p}.
\end{align*}

Thus, we have:
\begin{align*}
	\Vert (\Delta_{2} - \Delta_{1}) v  \Vert_{g_{2},n-1,p}
	&\leq
	c(n,p)
	\Vert (g_{2})^{-1} - (g_{1})^{-1} \Vert_{g_{2},n,p} \Vert (\nabla_{1})^2_{ij} v \Vert_{g_{2},n-1,p}
\\
	&\phantom{\leq}
	+ c(n,p)
	\Vert (g_{2})^{-1} \Vert_{g_{2},n,p} \Vert \triangle_{ij}^k \Vert_{g_{2},n-1,p} \Vert \ed v \Vert_{g_{2},n,p}.
\end{align*}

Therefore, if $\Vert g_{2} - g_{1} \Vert_{g_{1},n,p} \leq \delta < 1$, we obtain the following estimates:
\begin{align*}
	&
	\Vert (\Delta_{2} - \Delta_{1}) v  \Vert_{g_{2},n-1,p}
	\leq
	c(n,p)
	\delta \Vert v \Vert_{g_{2},n+1,p},
\\
	&
	\Vert v- \barv \Vert^{n+1,p}_{g_{2}}
	\leq
	c(n,p) \delta \Vert v \Vert_{g_{2},n+1,p},
\\
	&
	\Vert v- \barv \Vert^{n+1,p}_{g_{1}}
	\leq
	c(n,p) \delta \Vert v \Vert_{g_{1},n+1,p}.
\end{align*}

The lemma is thus proved.
\end{proof}

Suppose that $\{v_1, v_2, v_3\}$ is an orthonormal basis of the space $\Lambda_{k=1}^1$, and let $\{\barv_1, \barv_2, \barv_3\}$ denote the projections of $\{v_1, v_2, v_3\}$ onto the space $\Lambda_{k=1}^2$ with respect to the metric $g_{2}$. Consequently, $\{\barv_1, \barv_2, \barv_3\}$ forms an almost orthonormal basis for $\Lambda_{k=1}^2$.

\begin{lemma}\label{lem 4.2}
Suppose that $\Vert g_{2} - g_{1} \Vert_{g_{1},n,p} \leq \delta < 1$, where $n \geq n_p$. Let $\{ \barv_1, \barv_2, \barv_3\}$ be the projection of a unit orthogonal basis $\{v_1, v_2, v_3\}$ of $\Lambda_{k=1}^1$ onto $\Lambda_{k=1}^2$ with respect to the metric $g_{2}$. Then, there exists a constant $c_1(n,p)$ such that
\begin{align*}
	\big| \int_{\mathbb{S}^2} \barv_k \barv_l \, \dvol_{g_{2}} - \delta_{kl} \big| 
	\leq 
	c_1(n,p) \delta.
\end{align*}
Furthermore, there exists a unit orthogonal basis $\{\barv'_1, \barv'_2, \barv'_3\}$ of $\Lambda_{k=1}^2$ and a constant $c_2(n,p)$ such that
\begin{align*}
	\big| \barv_k - \barv'_k \big|, 
	\quad 
	\big| v_k - \barv'_k \big|, 
	\quad 
	\Vert \barv_k - \barv'_k \Vert_{g_{a},n+1,p}, 
	\quad 
	\Vert v_k - \barv'_k \Vert_{g_{a},n+1,p} \leq c_2(n,p) \delta,
	\quad a = 1, 2.
\end{align*}
\end{lemma}

\begin{proof}
The first inequality follows from Lemma \ref{lem 5.1r}, the Sobolev inequality, and the following formula:
\begin{align*}
	\int_{\mathbb{S}^2} \barv_k \barv_l \, \dvol_{g_{2}} 
	&
	= 
	\int_{\mathbb{S}^2} v_k v_l \, \dvol_{g_{1}} 
	+ 
	\int_{\mathbb{S}^2} v_k v_l (\dvol_{g_{2}} - \dvol_{g_{1}}) 
\\
	&
	\phantom{=\ } 
	+ \int_{\mathbb{S}^2} \big[ (\barv_k - v_k) \barv_l 
	+ v_k (\barv_l - v_l) \big] \, \dvol_{g_{2}}.
\end{align*}
The Gram-Schmidt process applied to $\{\barv_k\}_{k=1,2,3}$ yields the desired unit orthogonal basis. The second inequality follows from the almost orthonormality of $\{\barv_k\}_{k=1,2,3}$.
\end{proof}

Now, we prove that if two round metrics $g_{a}, a=1,2$ on $\Sigma$ are close, then their embeddings into the Euclidean sphere $(\mathbb{S}^2, \circg)$ are close to a rigid motion.

\begin{proposition}\label{pro 5.3r}
Suppose that $g_{a}, a=1,2$ are two round metrics with Gauss curvature 1 on the sphere $\mathbb{S}^2$. Let $\Phi_a$ be an isometric embedding of $(\Sigma, g_{a})$ into the Euclidean sphere $(\mathbb{S}^2, \circg)$. If $\Vert g_{2} - g_{1} \Vert_{g_{1},n,p} \leq \delta < 1$, where $n \geq n_p$, then there exist a constant $c(n,p)$ and a rigid motion $O \in O(3)$ such that
\begin{align*}
	\Vert \Phi_2 \circ \Phi_1^{-1} - O \Vert_{\circg,n+1,p} &\leq c(n,p) \delta, 
\\
	\Vert \Phi_2 - O \circ \Phi_1 \Vert_{g_{1},n+1,p} &\leq c(n,p) \delta.
\end{align*}
\end{proposition}

\begin{proof}
Note that two isometric embeddings of a round sphere with Gauss curvature 1 into the Euclidean sphere differ by a rigid motion. Therefore, it is sufficient to prove the inequality for two specific isometric embeddings. 

First, choose a unit orthogonal basis $\{v_k\}_{k=1,2,3}$ of $\Lambda_{k=1}^1$ and let $\{\barv'_k\}_{k=1,2,3}$ be the unit orthogonal basis of $\Lambda_{k=1}^2$ as in Lemma \ref{lem 4.2}. We then define two isometric embeddings from $(\Sigma, g_{a})$ to $(\mathbb{S}^2, \circg)$ using the above two unit orthogonal bases, i.e.,
\begin{align*}
	\Phi_1(x) &= \sqrt{\frac{4\pi}{3}}(v_1(x), v_2(x), v_3(x)), \\
	\Phi_2(x) &= \sqrt{\frac{4\pi}{3}}(\barv'_1(x), \barv'_2(x), \barv'_3(x)).
\end{align*}
We now consider the difference
\begin{align*}
	\Phi_2 - \Phi_1 = \sqrt{\frac{4\pi}{3}}(\barv'_1 - v_1, \barv'_2 - v_2, \barv'_3 - v_3).
\end{align*}
By Lemma \ref{lem 4.2}, we have
\begin{align*}
	\Vert (\barv'_k - v_k)\Vert_{g_1,n+1,p} \leq c(n,p) \delta, \quad k = 1, 2, 3.
\end{align*}
Thus, the proposition follows.
\end{proof}

\begin{remark}\label{rem 5.4r}
Proposition \ref{pro 5.3r} is analogous to Proposition 4.3 on almost isometries of $\mathbb{S}^2$ in \cite{KS2022}. The regularity assumption in Proposition 4.3 of \cite{KS2022} only requires an $L^\infty$ norm condition. In the proof of that proposition, the work \cite{FJM2002} is used. We show that by imposing a stronger regularity assumption in the present context, a more straightforward proof of Proposition \ref{pro 5.3r} can be derived.
\end{remark}

\section{Stability of semi-effective uniformization}\label{sec 6r}
We state the stability of semi-effective uniformizations of $2$-spheres as follows.

\begin{theorem}\label{thm 6.1r}
Let $p > 2$. Suppose that $(\Sigma, g_a)$, for $a = 1, 2$, are two $W^{2,p}$ spheres with Gauss curvature $K_a$. Assume that $(\Sigma, g_1)$ has the area $4\pi$ and two metrics are close in the following sense
\begin{align*}
    \Vert g_2 - g_1 \Vert_{g_1,2,p} 
    \leq 
    \delta,
\end{align*}
for some small $\delta > 0$.

Let $q \in \Sigma$ be a point on the sphere. Denote by $\Omega_a$ the uniformization conformal factor, where $\Omega_a(q) = 1$ and $\ed \Omega_a(q) = 0$. There exists a constant $\delta_0$ sufficiently small such that, if $\delta \leq \delta_0$, the following results hold.

Moreover, there exist constants $c_1(p)$, $c_2(p)$, and $c_3(p)$, which depend not only on $p$, but also on the following geometric quantities:
\begin{align*}
    \sup |K_1|, 
    \quad 
    \mathrm{diam}(\Sigma, g_1),
\end{align*}
or, if $k > 0$, on the quantities
\begin{align*}
    k_{1,m} = \min_{\mathbb{S}^2} K_1,
    \quad
    k_{1,M} = \max_{\mathbb{S}^2} K_1,
\end{align*}
such that
\begin{enumerate}[label=\alph*.]
\item
The uniformization conformal factor $\Omega_a$ are close in the sense that
\begin{align*}
	\big| \log\big(\frac{\Omega_2}{\Omega_1}\big) \big| \leq c(p) \delta,
	\quad
	\big\| \log\big(\frac{\Omega_2}{\Omega_1}\big) \big\|^{2,p}_{g_1} \leq c(p) \delta.
\end{align*}
The round metrics $(\Omega_{a})^{-2} \cdot g_{a}$ are close in the sense that
\begin{align*}
	\Vert (\Omega_{2})^{-2} \cdot g_{2} - (\Omega_{1})^{-2} \cdot g_{1} \Vert_{g_{1},2,p} 
	\leq 
	c_1(p) \delta
	< 1.
\end{align*}

\item
Let $\Phi_a$ be the isometric embedding of $(\Sigma, (\Omega_{a})^2 \cdot g_{a})$ into the Euclidean sphere $(\mathbb{S}^2, \circg)$. Then there exists a rigid motion $O \in O(3)$ such that
\begin{align*}
	\Vert \Phi_2 \circ \Phi_1^{-1} - O \Vert_{\circg,3,p} 
	\leq 
	c_2(p) \delta, 
	\quad
	\Vert \Phi_2 - O \circ \Phi_1 \Vert_{g_{1},3,p} 
	\leq 
	c_2(p) \delta.
\end{align*}

\end{enumerate}
\end{theorem}
\begin{proof}
\begin{enumerate}[label=\textit{\alph*.}]

\item
By the bound on $\Omega_1$ in Theorem \ref{thm 2.1r}, we have
\begin{align*}
	\Vert (\Omega_1)^{-2} \cdot g_2 - (\Omega_1)^{-2} \cdot g_1 \Vert_{(\Omega_1)^{-2} \cdot g_1,2,p}
	\leq 
	c(p) \delta.
\end{align*}
Thus, the Gauss curvature of $(\Sigma, (\Omega_1)^{-2} \cdot g_2)$ is close to 1, and we obtain the estimate
\begin{align*}
	\Vert K_{(\Omega_1)^{-2} \cdot g_2} - 1 \Vert_{(\Omega_1)^{-2} \cdot g_2,0,p}
	\leq
	c(p) \delta.
\end{align*}
Note that $\frac{\Omega_2}{\Omega_1}$ is the conformal factor for the metric $(\Sigma, (\Omega_1)^{-2} \cdot g_2)$, and the corresponding conformal transformation has constant Gauss curvature 1. We have the following conditions on the conformal factor:
\begin{align*}
	\frac{\Omega_2}{\Omega_1}(q) = 1, 
	\quad 
	\ed \big( \frac{\Omega_2}{\Omega_1} \big)(q) = 0.
\end{align*} 
Therefore, by Theorem \ref{thm 3.1r}, we deduce that
\begin{align*}
	\big| \log\big(\frac{\Omega_2}{\Omega_1}\big) \big| \leq c_0(p) \delta,
	\quad
	\big\| \log\big(\frac{\Omega_2}{\Omega_1}\big) \big\|^{2,p}_{g_1} \leq c_0(p) \delta.
\end{align*}
Hence, we conclude the following estimate:
\begin{align*}
	&\phantom{\leq}
	\Vert (\Omega_2)^{-2} \cdot g_2 - (\Omega_1)^{-2} \cdot g_1 \Vert_{g_1,2,p}
\\
	&
	\leq
	\Vert (\Omega_2)^{-2} \cdot g_2 - (\Omega_1)^{-2} \cdot g_2 \Vert_{g_1,2,p}
	+ \Vert (\Omega_1)^{-2} \cdot g_2 - (\Omega_1)^{-2} \cdot g_1 \Vert_{g_1,2,p}
\\
	&
	\leq
	c_1(p) \delta.
\end{align*}

\item
This is a direct consequence of Proposition \ref{pro 5.3r}.

\end{enumerate}
\end{proof}
\begin{remark}\label{rem 6.2r}
In comparison with the analogous Theorem 4.1 in \cite{KS2022}, the above theorem assumes lower regularity. Moreover, the above theorem is coordinate-independent, addressing a question raised in \cite{KS2022}.
\end{remark}

\section*{Acknowledgements}
The author is supported by the National Natural Science Foundation of China under Grant No. 12201338.

\appendix
\section{On Sobolev spaces $H_1^n(\mathbb{S}^2, \circg)$ and $\widehat{H}_2^n(\mathbb{S}^2, \circg)$}\label{appen Ar}

\begin{definition}\label{def A.1r}
Let $n \in \mathbb{N}$. $\pmb{H_1^n (\mathbb{S}^2, \circg)}$: the Hilbert space of $1$-forms on the unit sphere $(\mathbb{S}^2, \circg)$, equipped with the inner product
\begin{align*}
	\bdlangle \omega_1, \omega_2 \bdrangle_n
	=
	\sum_{k=0}^n 
	\int_{\mathbb{S}^2} 
	\langle \circnabla^k \omega_1, \circnabla^k \omega_2 \rangle_{\circg} \, \dvol_{\circg}
\end{align*}
$\pmb{\widehat{H}_2^n (\mathbb{S}^2, \circg)}$: the Hilbert space of trace-free symmetric $2$-tensors on $(\mathbb{S}^2, \circg)$, with the inner product
\begin{align*}
	\bdlangle \widehat{\Xi}_1, \widehat{\Xi}_2 \bdrangle_n
	=
	\sum_{k=0}^n 
	\int_{\mathbb{S}^2} 
	\langle \circnabla^k \widehat{\Xi}_1, \circnabla^k \widehat{\Xi}_2 \rangle_{\circg} \, \dvol_{\circg}
\end{align*}
In the following, we will use the shorthand notations $H_1^n$ and $\widehat{H}_2^n$ to denote these spaces.
\end{definition}

Let $\{Y_l^m\}_{l \in \mathbb{N}, -l \leq m \leq l }$ denote the set of spherical harmonic functions of degree $l$ and order $m$, which forms an orthonormal basis of the $L^2$ function space on the unit sphere. Using this, we can construct orthonormal bases for $H_1^n$ and $\widehat{H}_2^n$.

\begin{proposition}\label{prop A.2r}
$\{ \ed Y_l^m, \star\ed Y_l^m \}_{l \in \mathbb{Z}_{+}, -l \leq m \leq l }$ is an orthonormal basis of $H_1^n$, and
\begin{align*}
	\bdlangle \ed Y_l^m, \ed Y_l^m \bdrangle_n
	\sim
	\lambda_l^{n+1}
	\sim
	l^{2(n+1)},
	\quad
	\bdlangle \star \ed Y_l^m, \star \ed Y_l^m \bdrangle_n
	\sim
	\lambda_l^{n+1}
	\sim
	l^{2(n+1)}.
\end{align*}
$\{ \calL \ed Y_l^m \}_{l \in \mathbb{Z}_{\geq 2}, -l \leq m \leq l } \cup \{ \calL \star\ed Y_l^m \}_{l \in \mathbb{Z}_{+}, -l \leq m \leq l }$ is an orthonormal basis of $\widehat{H}_2^n$, and
\begin{align*}
	\bdlangle \calL \ed Y_l^m, \calL \ed Y_l^m \bdrangle_n
	\sim
	\lambda_l^{n+2}
	\sim
	l^{2(n+2)},
	\quad
	\bdlangle \calL \star \ed Y_l^m, \calL \star \ed Y_l^m \bdrangle_n
	\sim
	\lambda_l^{n+2}
	\sim
	l^{2(n+2)}.
\end{align*}
\end{proposition}

We now introduce the dual spaces of $H_1^n$ and $\widehat{H}_2^n$ with respect to the $L^2$ inner product.
\begin{definition}\label{def A.3r}
Let $n\in \mathbb{N}$. $\pmb{H_1^{-n} (\mathbb{S}^2, \circg)}$: the dual Hilbert space of $H_1^n$ with respect to $L^2$ inner product, equipped with the norm
\begin{align*}
	\| \omega \|_{\circg,-n,2}
	=
	\sup_{\omega' \neq 0} \frac{\bdlangle \omega, \omega' \bdrangle_0}{\| \omega' \|_{\circg,n,2}}
\end{align*}
$\pmb{\widehat{H}_2^{-n} (\mathbb{S}^2, \circg)}$: the dual Hilbert space of $\widehat{H}_2^n$ with respect to $L^2$ inner product, with the norm
\begin{align*}
	\| \widehat{\Xi} \|_{\circg,-n,2}
	=
	\sup_{\widehat{\Xi}' \neq 0} \frac{\bdlangle \widehat{\Xi}, \widehat{\Xi}' \bdrangle_0}{\| \widehat{\Xi}' \|_{\circg,n,2}}
\end{align*}
\end{definition}
We also have a similar proposition regarding the bases for $H_1^{-n}$ and $\widehat{H}_2^{-n}$.
\begin{proposition}\label{prop A.4r}
$\{ \ed Y_l^m, \star\ed Y_l^m \}_{l \in \mathbb{Z}_{+}, -l \leq m \leq l }$ is an orthonormal basis of $H_1^{-n}$, and
\begin{align*}
	\bdlangle \ed Y_l^m, \ed Y_l^m \bdrangle_{-n}
	\sim
	\lambda_l^{1-n}
	\sim
	l^{2-2n},
	\quad
	\bdlangle \star \ed Y_l^m, \star \ed Y_l^m \bdrangle_n
	\sim
	\lambda_l^{1-n}
	\sim
	l^{2-2n}.
\end{align*}
$\{ \calL \ed Y_l^m \}_{l \in \mathbb{Z}_{\geq 2}, -l \leq m \leq l } \cup \{ \calL \star\ed Y_l^m \}_{l \in \mathbb{Z}_{+}, -l \leq m \leq l }$ is an orthonormal basis of $\widehat{H}_2^n$, and
\begin{align*}
	\bdlangle \calL \ed Y_l^m, \calL \ed Y_l^m \bdrangle_{-n}
	\sim
	\lambda_l^{2-n}
	\sim
	l^{2(2-n)},
	\quad
	\bdlangle \calL \star \ed Y_l^m, \calL \star \ed Y_l^m \bdrangle_{-n}
	\sim
	\lambda_l^{2-n}
	\sim
	l^{2(2-n)}.
\end{align*}
\end{proposition}

\Address


\begin{thebibliography}{100}

\bibitem[CY87]{CY1987}
Chang, S.-Y. A.; Yang, P.:
Prescribing Gaussian curvature on $\mathbb{S}^2$.
\textit{Acta Math.} \textbf{159} (1987), 215–259.

\bibitem[CK93]{CK1993}
Christodoulou, D.; Klainerman, S.:
\emph{The global nonlinear stability of the Minkowski space.}
Princeton Mathematical Series, 41. \emph{Princeton University Press, Princeton, NJ}, 1993. x+514 pp.

\bibitem[FJM02]{FJM2002}
Friesecke, G.; James, R., M\"uller,S.:
A theorem on geometric rigidity and the derivation of nonlinear plate theory from three-dimensional elasticity.
\emph{Comm. Pure Appl. Math.} \textbf{55} (2002), 1461–1506.

\bibitem[J61]{J1961}
John, F.:
Rotation and strain.
\emph{Commun.Pure Appl. Math.} \textbf{44} (1961), 391–413.

\bibitem[KS22]{KS2022}
Klainerman, S.; Szeftel, J.:
Effective results on uniformization and intrinsic GCM spheres in perturbations of Kerr. With an appendix by Camillo De Lellis. 
\emph{Ann. PDE} \textbf{8} (2022), no. 2, Paper No. 18, 89 pp.

\bibitem[L18]{L2018}
Le, P.:
\emph{The perturbation theory of null hypersurfaces and the weak null Penrose inequality},
DISS. ETH Nr. 25387. \url{https://doi.org/10.3929/ethz-b-000334917}

\bibitem[O82]{O1982}
Onofri, E.:
On the positivity of the effective action in a theory of random surfaces. 
\textit{Comm. Math. Phys.} \textbf{86} (1982), no. 3, 321–326. 

\bibitem[SWW09]{SWW2009}
Shi, Y.; Wang, G.; Wu, J.:
On the behavior of quasi-local mass at the infinity along nearly round surfaces.
\emph{Ann. Global Anal. Geom.} \textbf{36} (2009), no.4, 419–441.

\bibitem[W23]{W2023}
Wolff, M.:
A De Lellis–M\"uller type estimate on the Minkowski lightcone.
\textit{Calc. Var. Partial Differential Equations} \textbf{63} (2024), no. 7, Paper No. 185, 30 pp.

\end{thebibliography}
\end{document}